\documentclass{article}
\usepackage[dvipsnames,svgnames,x11names]{xcolor}
\usepackage{tikz,geometry,amsmath,hyperref,bm,flowchart,booktabs,overpic}
\usepackage[utf8]{inputenc}
\usepackage{algorithm,algorithmic}
\usepackage{authblk}
\usetikzlibrary{shapes.geometric,arrows}
\tikzset{>=stealth}
\hypersetup{hidelinks}

\newcommand\keywords[1]{\small{\textbf{Keywords}: #1}}

\title{On solution of conformal mapping for a lower half plane containing a symmetrical noncircular cavity}
\author[1]{Luobin Lin}
\author[2]{Fuquan Chen\thanks{Corresponding author: phdchen@fzu.edu.cn (Fuquan Chen)}}
\author[1]{Xianhai Huang}
\date{}
\affil[1]{Fujian Provincial Key Laboratory of Advanced Technology and Informatization in Civil Engineering, College of Civil Engineering, Fujian University of Techonology, No. 69 Xueyuan Road, Shangjie University Town, Fuzhou, 350118, Fujian, China}
\affil[2]{College of Civil Engineering, Fuzhou University, No. 2 Xueyuan Road, Shangjie University Town, Fuzhou, 350108, Fujian, China}

\begin{document}
\maketitle

\begin{abstract}
In this paper, we provide a candidate solution to obtain the coefficients of the conformal mapping for a lower half plane containing a symmetrical noncircular cavity using penalty function method and modified Particle Swarm Method. The nonconvexity of the penalty function is proven via the concept of convex function and proof of contradiction. The solution procedure is presented very detailedly in pseudocodes to ensure that the solution can be fully repeated and further improved. The solution accuracy and efficiency are also discussed.
\end{abstract}


\section{Introduction}
\label{sec.1}

Analytic methods are widely used in mechanical analyses in tunnel engineering. However, owing to the complicated shapes of tunnels, it is difficult to analyze mechanical behaviour of tunnel and surrounding geomaterial in a physical plane. Thus, mechanical analyses generally focused on circular tunnels in early studies \cite{Mindlin1940,Green1940,LingChinBing1948,Lekhnitskii1963,Sagaseta_image_method}. Conformal mapping is introduced to mechanical analyses in tunnel engineering to transform a physical plane with a complicated boundary onto a mapping plane with a circular boundary, which greatly facilites the analyses of mechanical behaviour of noncircular tunnels. 

With the aid of conformal mapping, the complex variable method \cite{Muskhelishvili1966} is systematically proposed as a very flexible analytic method in 2-dimensional mechanical analyses in tunnel engineering, and is suitable to various shapes of tunnel boundaries, such as rectangular \cite{Huo2006}, notched circular \cite{Exadaktylos2002}, and U-shaped \cite{Lu2015,Lu2014}. The complex variable method can also be used for multi-tunnel problems together with conformal mappings \cite[]{Lu1997,Self2019deep_twin_liner,fanghuangcheng2021}. An alternative conformal mapping function \cite[]{Howland1939} is proposed and modified for noncircular twin tunnels at great depth \cite[]{Zeng2018}. However, the above researches mainly focus on deep tunnels, owing to the lack of pathbreaking conformal mapping functions for shallow tunnels, which can be abstracted as a lower half plane containing a cavity.

In 1997, Verrujit \cite{Verruijt1997traction,Verruijt1997displacement} proposed a compact conformal mapping function that maps a lower half plane containing a circular cavity into a unit ring. Such a conformal mapping function greatly extends the usage of complex variable method and inspires many researches \cite{Wanglizhong2009,zhang2018complex,Self2021APM_shallow_twin,Wanghuaning2021shallow_rigid_lining}. Furthermore, the gravity gradient and the unbalanced resultant due to tunnel excavation can be considered with the Verruijt conformal mapping function \cite[]{Strack2002phdthesis,Strack_Verruijt2002buoyancy,Verruijt_Strack2008buoyancy,Lu2016,Lu2019new_solution,Self2020JEM}.

All those researches of shallow tunnels still focus on circular cavities, instead of more commonseen noncircular cavities in tunnel engineering, because no conformal mapping function for a lower half plane containing a cavity of an arbitrary shape is proposed.  Gratifyingly, Zeng et al. \cite{zengguisen2018thesisenglish,Zengguisen2019} propose an elegant conformal mapping function for a lower half plane containing a noncircular cavity, which further extends the usage of complex variable method. Unfortunately, the solution method and procedure of this conformal mapping function is described very briefly, except that the optimization theory is applied. Referring to Ref \cite[]{Zengguisen2019}, it is difficult to know 1) whether or not the objective funtion is convex or nonconvex; 2) via which method the conformal mapping function is solved; 3) the solution accuracy and efficiency.

The lack of solution procedure increases the cost to understand this conformal mapping function, and greatly limits its usage. This paper is focusing on providing a candidate solution method and corresponding solution procedure in detail, so that the solution of this conformal mapping function can be repeated and more widely used in complex variable method.

\section{Conformal mapping expression and problem definition}
\label{sec.2}

\begin{figure}[hbt]
  \centering
  \begin{tikzpicture}
    \tikzstyle{every node}=[scale=0.7]
    \shade[top color=gray!10, bottom color=gray!40](-8,0)rectangle(-1,-6);
    \draw[Emerald,line width=1pt](-8,0)--(-1,0);
    \draw[dash dot, red](-4.5,0) -- (-4.5,-6);
    \fill[white](-4.5,-2.7)ellipse[x radius=2.5,y radius=2];
    \draw[red,line width=1pt](-4.5,-2.7)ellipse[x radius=2.5,y radius=2];
    \draw[brown, dashed , line width=2pt] (-4.5,-4.7) arc [x radius=2.5, y radius=2, start angle=-90, end angle=90];
    \draw[black, dashed, line width=1pt] (-4.5,-2.7) circle [radius=2];
    \draw[black, ->, line width=1pt] (-4.5,-2.7) -- (-2.5,-2.7) node [above left] {$R_{\min}$};
    \draw[black, dashed] (-4.5,-2.7) -- (-5.5,-2.7);
    \draw[black, <->] (-5,0) -- (-5,-2.7);
    \node at (-5,-1.35) [rotate=90,above] {$H_{\rm mid}$};
    \filldraw[black] (-4.5,-2.7) circle [radius=0.05];
    \draw[->](-4.5,0)--(-3.5,0)node[above]{$x$};
    \draw[->](-4.5,0)--(-4.5,1)node[right]{$y$};
    \node at(-4.5,0)[below]{$O$}; \filldraw [blue] (-4.5,-4.7) circle [radius=0.05];
    \node at (-4.5,-4.7) [above,blue] {$P_0$};
    \filldraw [blue] (-4.5,-0.7) circle [radius=0.05];
    \node at (-4.5,-0.7) [below,blue] {$P_\pi$};
    \node at(-4.5,-4.7)[below,red]{Tunnel boundary};
    \node at(-4.5,0)[above left,Emerald]{Ground surface};
    \node at(-4.5,-6)[below]{a) Physical plane};

    \fill[outer color=gray!10, inner color=gray!40](5.5,-3)circle[radius=3];
    \draw[Emerald,line width=1pt](5.5,-3)circle[radius=3];
    \draw[dash dot, red] (2.5,-3) -- (8.5,-3); \fill[white](5.5,-3)circle[radius=1];
    \draw[red,line width=1pt](5.5,-3)circle[radius=1];
    \draw[brown, dashed, line width=2pt] (6.5,-3) arc [radius=1, start angle =0, end angle=180];
    \draw[->](5.5,-3)--(6.3,-3)node[above]{$\rho$};
    \draw[->](5.8,-3)arc[start angle=0,end angle=360,radius=0.3];
    \node at(5.8,-3)[above right]{$\theta$};
    \node at(5.5,-3)[below]{$o$}; \filldraw [blue] (6.5,-3) circle [radius=0.05];
    \node at (6.5,-3) [right,blue] {$P_0$};
    \filldraw [blue] (4.5,-3) circle [radius=0.05];
    \node at (4.5,-3) [left,blue] {$P_\pi$};
    \draw[->,Emerald](5.5,-3)--({5.5+3*cos(135)},{-3+3*sin(135)})node[right]{\textcolor{Emerald}{1}};
    \draw[->,red](5.5,-3)--({5.5+cos(-65)},{-3+sin(-65)})node[below right] {\textcolor{red}{$\alpha$}};
    \node at(5.5,-6)[below]{b) Mapping plane};

    \draw[<->,line width=2pt](-0.8,-3)--(2.2,-3);
    \node at(0.7,-3)[above,align=left]{\parbox{10cm}{
        \begin{equation*}
          \begin{aligned}
            \omega(\zeta)=
            & -{\rm i}c_{0} \frac{1+\zeta}{1-\zeta}                               \\
            & +{\rm i}\sum\limits_{k=1}^{m}c_k\left(\zeta^k-\zeta^{-k}\right)
          \end{aligned}
        \end{equation*}
      } }; \node at(1.3,-0.5){\parbox{10cm}{
        \begin{equation*}
          \left\{
            \begin{aligned}
              \textcolor{Emerald}{x} & =c_{0} \frac{\sin\textcolor{Emerald}{\theta}}{1-\cos\textcolor{Emerald}{\theta}}-\sum\limits_{k=1}^{m}c_k\sin{k}\textcolor{Emerald}{\theta} \\
              \textcolor{Emerald}{y} & =0
            \end{aligned}
          \right.
        \end{equation*}
      } }; \node at(0.9,-4.8){\parbox{10cm}{
        \begin{equation*}
          \left\{
            \begin{aligned}
              \textcolor{red}{x}=
              & c_{0} \frac{2\textcolor{red}{\alpha}\sin\textcolor{red}{\theta}}{1+\textcolor{red}{\alpha}^2-2\textcolor{red}{\alpha}\cos\textcolor{red}{\theta}} \\
              & -\sum\limits_{k=1}^{m}{c_k}\left(\textcolor{red}{\alpha}^k+\textcolor{red}{\alpha}^{-k}\right)\sin{k}\textcolor{red}{\theta}                     \\
              \textcolor{red}{y}=
              & -c_{0} \frac{1-\textcolor{red}{\alpha}^2}{1+\textcolor{red}{\alpha}^2-2\textcolor{red}{\alpha}\cos\textcolor{red}{\theta}}                           \\
              & +\sum\limits_{k=1}^{m}{c_k}\left(\textcolor{red}{\alpha}^k-\textcolor{red}{\alpha}^{-k}\right)\cos{k}\textcolor{red}{\theta}
            \end{aligned}
          \right.
        \end{equation*}
      } };
  \end{tikzpicture}
  \caption{Conformal mapping of the elastic region into a unit ring}
  \label{fig.1}
\end{figure}
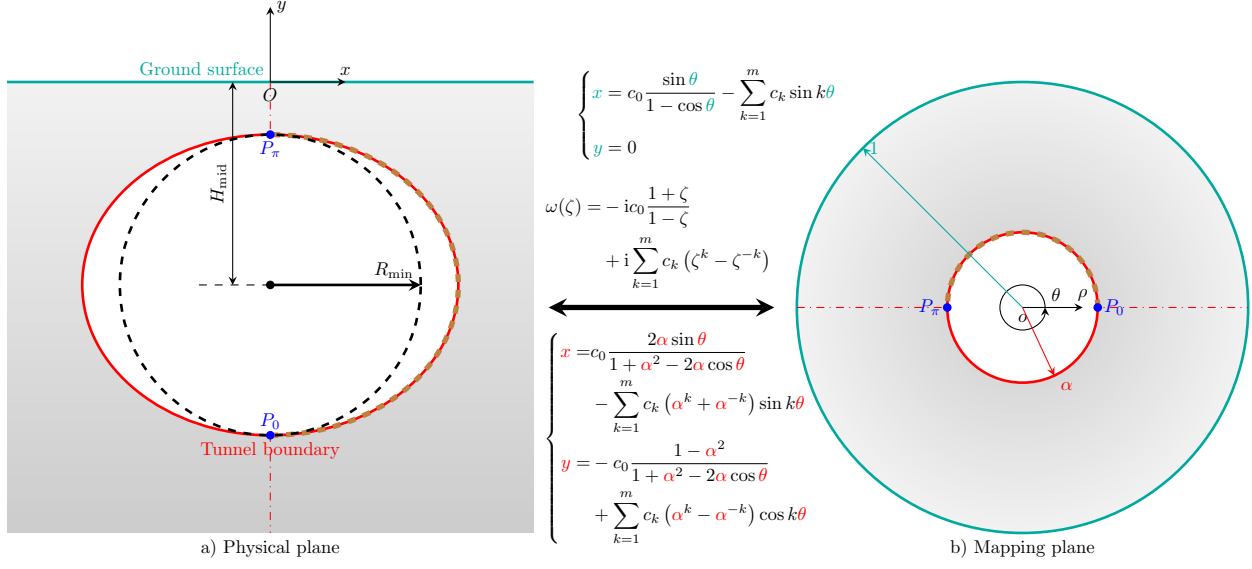

The conformal mapping function proposed by Zeng et al. \cite{zengguisen2018thesisenglish, Zengguisen2019} is as follows:
\begin{equation}
  z=\omega(\zeta)={\rm i}c_{0} \frac{1+\zeta}{1-\zeta}+{\rm i}\sum\limits_{k=1}^{n}c_{k}(\zeta^{k}-\zeta^{-k})
  \label{eq.1}
\end{equation}

\noindent where $z$ denotes complex coordinate in the physical plane, $\zeta$ denotes complex coordinate in the mapping plane, ${\rm i}$ denotes imaginary unit, $c_k (k=0,1,2,3,\cdots,n)$ denote real coefficients to be determined. As shown in Fig. \ref{fig.1}, via Eq. (\ref{eq.1}), the elastic region in the rectangular coordinate system $xOy$ (physical plane) is mapped onto a unit ring with an inner radius of $\alpha$ in the polar coordinate system $\rho o \theta$ (mapping plane) with $\zeta=\rho\cdot\sigma=\rho\cdot{\rm e}^{{\rm i}\theta}$. The conformal mapping function in Eq. (\ref{eq.1}) can be rewritten in the polar coordinate form:

\begin{equation}
  \omega(\rho,\theta)={\rm i}c_{0} \frac{1+\rho {\rm e}^{{\rm i}\theta}}{1-\rho {\rm e}^{{\rm i}\theta}}+{\rm i}\sum\limits_{k=1}^{n}c_{k}(\rho^k {\rm e}^{{\rm i}k\theta}-\rho^{-k} {\rm e}^{-{\rm i}k\theta})
  \label{eq.2}
\end{equation}

\noindent and in the rectangular coordinate form:

\begin{equation}
  \left\{
    \begin{aligned}
      & x(\rho,\theta)=-c_{0} \frac{2\rho\sin\theta}{1+\rho^2-2\rho\cos\theta}-\sum\limits_{k=1}^{n}{c_{k}}\left(\rho^k+\rho^{-k}\right)\sin{k}\theta \\
      & y(\rho,\theta)=c_{0} \frac{1-\rho^2}{1+\rho^2-2\rho\cos\theta}+\sum\limits_{k=1}^{n}{c_{k}}\left(\rho^k-\rho^{-k}\right)\cos{k}\theta
    \end{aligned}
  \right.
  \label{eq.3}
\end{equation}

When $\rho=\alpha$, the conformal mapping of the elasto-plastic boundary in the polar coordinate form and in the rectangular coordinate form can be obtained, respectively:

\begin{equation}
  \omega(\alpha,\theta)={\rm i}c_{0} \frac{1+\alpha {\rm e}^{{\rm i}\theta}}{1-\alpha {\rm e}^{{\rm i}\theta}}+{\rm i}\sum\limits_{k=1}^{n}c_{k}(\alpha^k {\rm e}^{{\rm i}k\theta}-\alpha^{-k} {\rm e}^{-{\rm i}k\theta})
  \label{eq.4}
\end{equation}

\begin{equation}
  \left\{
    \begin{aligned}
      & x(\alpha,\theta)=-c_{0} \frac{2\alpha\sin\theta}{1+\alpha^2-2\alpha\cos\theta}-\sum\limits_{k=1}^{n}{c_k}\left(\alpha^k+\alpha^{-k}\right)\sin{k}\theta \\
      & y(\alpha,\theta)=c_{0} \frac{1-\alpha^2}{1+\alpha^2-2\alpha\cos\theta}+\sum\limits_{k=1}^{n}{c_k}\left(\alpha^k-\alpha^{-k}\right)\cos{k}\theta
    \end{aligned}
  \right.
  \label{eq.5}
\end{equation}

According to the description of the solution procedure in Ref \cite{Zengguisen2019}, a series of sample points along the right half of the inner boundary of the elastic region in the physical plane are selected:

\begin{equation}
  \label{eq.6}
  \left\{
    \begin{aligned}
      & X_{i} = -c_{0} \frac{2\alpha\sin\theta_{i}}{1+\alpha^2-2\alpha\cos\theta_{i}}-\sum\limits_{k=1}^{n}{c_k}\left(\alpha^k+\alpha^{-k}\right)\sin{k}\theta_{i} = x^{0}_{i} \\
      & Y_{i} = c_{0} \frac{1-\alpha^2}{1+\alpha^2-2\alpha\cos\theta_{i}}+\sum\limits_{k=1}^{n}{c_k}\left(\alpha^k-\alpha^{-k}\right)\cos{k}\theta_{i} = y^{0}_{i}\\
    \end{aligned}
  \right.
  ,\quad
  \begin{aligned}
    & i=1,2,3,\cdots,m \\
    & m \geq n+2
  \end{aligned}
\end{equation}

\noindent where $x^{0}_{i}$ and $y^{0}_{i}$ denote the horizontal and vertical coordinates of the sample points in the physical coordinate system $xOy$ (which are known), $\theta_{i}$ denote the polar angles in the mapping plane (which are unknown). When $\theta=0$ or $\pi$, Eqs. (\ref{eq.3}) and (\ref{eq.5}) turn to $x=0$, which are corresponding to points $i=1$ and $i=m$ in Eq. (\ref{eq.6}). Such a feature would not change for different $c_k$, and these two points are marked as $P_0$ and $P_\pi$ in the physical and mapping planes in Fig \ref{fig.1}, respectively. Thus, Eq. (\ref{eq.6}) is an over-determined nonlinear equation system containing $m+n$ variables and $2m-2$ nonlinear equations.

To solve these $m+n$ variables in Eq. (\ref{eq.6}) with the $2m-2$ nonlinear equations, the penalty function method of optimization is applied in Ref \cite{Zengguisen2019}. However, Ref \cite{Zengguisen2019} does not provide detailed penalty function expression, here the authors complete the expression according to the definition of penalty function:

\begin{equation}
  \label{eq.7}
  \begin{aligned}
    & \min f(x) \\
    & {\rm s.t.}
      \left\{
      \begin{aligned}
        & X_{i} - x_{i}^{0} = 0,\quad i=2,3,\cdots,m-1 \\
        & Y_{i} - y_{i}^{0} = 0,\quad i=1,2,3,\cdots,m \\
        & 0 < \alpha < 1 \\
        & \theta_{i} < \theta_{i+1},\quad i=2,3,4,\cdots,m-2 \\
        & \theta_{1} = 0 \\
        & \theta_{m} = \pi \\
      \end{aligned}
      \right.
  \end{aligned}
\end{equation}

\noindent where $f(x)$ denotes the objective function, which is not set according to Ref \cite{Zengguisen2019} ($f(x)=0$), because there is no objective function, and our aim is to solve Eq. (\ref{eq.6}). Thus, the following penalty function can be established for optimization:

\begin{equation}
  \label{eq.8}
  \begin{aligned}
    F({\bm C},{\bm \theta})=\varOmega\cdot & \left\{  \sum\limits_{i=2}^{m-1} \left[X_{i}({\bm C},{\bm \theta})-x_{i}^{0} \right]^{2} + \sum\limits_{i=1}^{m} \left[Y_{i}({\bm C},{\bm \theta}) - y_{i}^{0} \right]^{2}\right. \\
    & \left. + \left[\max\left(0, \alpha-1 \right)\right]^{2} + \left[\max(0, -\alpha)\right]^{2} + \sum\limits_{i=1}^{m-1}\left[\max (0,\theta_{i} - \theta_{i+1})\right]^{2} \right\} \\
  \end{aligned}
\end{equation}

\noindent where $\varOmega$ denotes penalty factor. To start the optimization, the following initial values are recommended in Ref \cite{Zengguisen2019}:

\begin{equation}
  \label{eq.9}
  \left\{
    \begin{aligned}
      & {\bm C}[1] = \alpha = \frac{\frac{D}{2}}{H+\sqrt{H^{2}-\left(\frac{D}{2}\right)^{2}}} \\
      & {\bm C}[2] = c_{0} = -H\cdot\frac{1-\alpha^{2}}{1+\alpha^{2}} \\
      & {\bm C}[i] = 0,\quad i=3,4,5,\cdots,n+2 \\
      & {\bm \theta}[i] = \pi\cdot\frac{i-1}{m-1},\quad i=1,2,3,\cdots,m
    \end{aligned}
  \right.
\end{equation}

\noindent where $D$ denotes the net height of the cavity, $H$ denotes the buried depth of the cavity.

Though Ref \cite{Zengguisen2019} presents the conformal mapping functions of several commonseen tunnel shapes in the lower half plane, the solution method and procedure are vague, which limits the usage of the elegant conformal mapping function. For this a reason, this paper is dedicated to discuss a candidate solution method and corresponding procedure of the optimization of the penalty function in Eq. (\ref{eq.8}).

Clearly, the form of the penalty function is complicated. To apply the optimization theory, the following questions should be answered:

\noindent(1) Is the penalty function in Eq. (\ref{eq.8}) convex or nonconvex?

\noindent(2) If the penalty function is convex, which local optimization method would be efficient?

\noindent(3) If the penalty function is nonconvex, how to optimize the nonconvex problem with a global optimization method?

The following text would answer the above questions.

\section{Verification of nonconvexity}
\label{sec.3}

By observing the form of the penalty function in Eq. (\ref{eq.8}), it is highly possible that the penalty function is nonconvex. Besides, owing to the complicated form of the penalty function, the gradient or Hessian matrix method would not be suitable to determine whether or not the penalty function is convex. For this reason, the definition of convex function and proof of contradiction are together applied to determine the convexity of the penalty fucntion with the aid of a program coded by the authors.

Assumption that the penalty function in Eq. (\ref{eq.8}) is convex is first made. According to the definition of convex function, the following inequality should always be satisfied for arbitrary ${\bm C}$ or ${\bm \theta}$:

\begin{equation}
  \label{eq.10}
  G = \lambda \cdot F({\bm C}[:,1],{\bm \theta}[:,1]) + (1-\lambda) \cdot F({\bm C}[:,2],{\bm \theta}[:,2]) - F({\bm C}[:,3],{\bm \theta}[:,3]) \geq 0
\end{equation}

\noindent where $\lambda$ denotes an arbitrary real numeric in the range $(0,1)$, ${\bm C}$ and ${\bm \theta}$ are expanded into matrices of sizes of $(n+2)\times 3$ and $m \times 3$, respectively, and

\begin{equation}
  \label{eq.11}
  \left\{
    \begin{aligned}
      & {\bm C}[:,3] = \lambda \cdot {\bm C}[:,1] + (1-\lambda) \cdot {\bm C}[:,2] \\
      & {\bm \theta}[:,3] = \lambda \cdot {\bm \theta}[:,1] + (1-\lambda) \cdot {\bm \theta}[:,2] \\
    \end{aligned}
  \right.
\end{equation}

\noindent ${\bm C}[:,1] $ and ${\bm C}[:,2] $, ${\bm \theta}[:,1] $ and ${\bm \theta}[:,2] $ denote two arbitrary sets of variables. As long as one counter-example for Eq. (\ref{eq.10}) can be found ($G<0$), the contradiction for the assumption that Eq. (\ref{eq.8}) is convex would be met, and the nonconvexity of the penalty function is proven.

Ref \cite{zengguisen2018thesisenglish} indicates that an elliptical cavity would be relatively simple for conformal mapping, and the following elliptical cavity would be applied for further contradiction proof:

\begin{equation}
  \label{eq.12}
  {\left(\frac{x}{A}\right)}^{2}+{\left(\frac{y+H}{B}\right)}^{2}=1
\end{equation}

\noindent where $A=3$, $B=2$, $H=5$. Thus, Eq. (\ref{eq.12}) is a specific elliptical curve, and the following sample points along this elliptical curve can be selected 

\begin{equation}
  \label{eq.13}
  \left\{
    \begin{aligned}
      & x^{0}_{i} = A \cdot \cos \phi_{i} \\
      & y^{0}_{i} = B \cdot \sin \phi_{i} - H \\
    \end{aligned}
  \right.
  ,\quad
  \begin{aligned}
    & \phi_{i}=\pi\cdot\frac{i}{m}-\frac{\pi}{2} \\
    & i=1,2,3,\cdots,m
  \end{aligned}
\end{equation}

\noindent where $x_{i}^{0}$ and $y_{i}^{0}$ denote horizontal and vertical coordinates of the sample points along the elliptical boundary of the cavity, respectively, $\phi_{i}$ denote the polar angle in the physical plane in the lower half plane, and $n=10$, $m=90$ would be accurate enough.

The pseudocodes of the verification procedure can be found in \texttt{Algorithm \ref{alg.1}}. In Lines 10, 11, 13 and 17, $t$ denotes a random numeric within the range $(0,1)$ generated by computer. In Lines 22 and 23, $\lambda$ is identical to that in Eqs. (\ref{eq.10}) and (\ref{eq.11}), and is another random numeric within the range $(0,1)$ generated by computer. In Line 11, the value $20$ can be replaced by arbitrary real number. In Line 13, the range $(-1,1)$ is based on the solutions of Tables 1-3 in Ref \cite{Zengguisen2019} and the solutions in Fig.3-2 and Table 6-1 in Ref \cite{zengguisen2018thesisenglish}, as well as the fact that the Laurent series in Eq. (\ref{eq.1}) would converge.

Substituting $A$, $B$, $H$, and $n$, $m$, $\varOmega=10^{10}$ into the pseudocodes in \texttt{Algorithm \ref{alg.1}} yields the iteration reps $N$ and the first $G$ that meets the demand $G<0$. Generally, the computation in \texttt{Algorithm \ref{alg.1}} is very fast (less than 0.1 second). In this paper, the pseudocodes are realized by \texttt{FORTRAN}. The existence of $G$ ($G<0$) violates Eq. (\ref{eq.10}) and indicates that the penalty function in Eq. (\ref{eq.8}) is nonconvex.

\section{Solution method and procedure}
\label{sec.4}

Since the penalty function in Eq. (\ref{eq.8}) has been proven nonconvex, we should find an efficient method to gloablly minimize $F({\bm C},{\bm \theta})$ to approach zero and to obtain the corresponding solution ${\bm C}^{*}$ and ${\bm \theta}^{*}$.  To be consistent with previous sections, the to-be-calculated example is also chosen as the elliptical curve in Eqs. (\ref{eq.12}) and (\ref{eq.13}).

This paper aims at presenting the solution method and procedure of the optimization for the problem defined in Section \ref{sec.2}, instead of presenting method comparison or competetion. Therefore, only one global would be used and detailedly discussed to serve as a modest spur. The solution method and procedure should be concretely and detailedly presented via specific calculation example, and the pseudocodes and corresponding real code should be presented as well, thus, the solution method and procedure can be repeated and possibly further improved.

Among many global optimization methods, the Particle Swarm Method \cite{kennedy1995particle} is applied. The screen strategy of the Particle Swarm Method always records global optimum variable set, as well as individual optimum variable sets, in one single iteration. Such feature would reduce the possibility of missing the correct global optimum. The original Particle Swarm Method is slightly modified to ensure convergence of the penalty function in Eq. (\ref{eq.8}). The pseudocodes of the modified Particle Swarm Method are presented in \textbf{Algorithm} \ref{alg.2}, in which the variable initialization, parameter updating strategies, and variable updating strategies can be found in \textbf{Algorithm} \ref{alg.3}-\ref{alg.5}, respectively. The symbols in \textbf{Algorithm} \ref{alg.2}-\ref{alg.5} are explained in the following text along with the usage of relavent pseudocodes.

\subsection{Algorithm explanation}
\label{sec.4.1}

The \textbf{Input} parameters in \textbf{Algorithm} \ref{alg.2} are explained as follows: $VTR$ denotes value to reach; $N_{\max}$ denotes the maximum iteration reps, which should not be large, since the method should be efficient; $S$ denotes the population size. $\varOmega$ is identical to that large real number in Eq. (\ref{eq.8}); ${\bm C}_{v\min}^{\rm init}$, ${\bm C}_{v\min}^{\rm end}$, ${\bm C}_{v\max}^{\rm init}$, ${\bm C}_{v\max}^{\rm end}$ denote the initial and ending limits for the minimum and maximum search velocities of ${\bm C}$, respectively; while ${\bm \theta}_{v\min}^{\rm init}$, ${\bm \theta}_{v\min}^{\rm end}$, ${\bm \theta}_{v\max}^{\rm init}$, ${\bm \theta}_{v\max}^{\rm end}$ are similarly defined for  ${\bm \theta}$, respectively. $w_{\rm init}$ and $w_{\rm end}$ denote the initial and ending inertias in each iteration, respectively; $d_{\rm 1init}$ and $d_{\rm 1end}$ denote the initial and ending individual learning factors, respectively; $d_{\rm 2init}$ and $d_{\rm 2end}$ denote the initial and ending global learning factors, respectively. ${\bm C}^{*}$ and ${\bm \theta}^{*}$ denote the results of the whole optimization procedure.

In \textbf{Algorithm} \ref{alg.2}-\ref{alg.5}, ${\bm C}$ and ${\bm \theta}$ are expanded into matrices of sizes of $(n+2)\times S$, and $m\times S$, respectively; ${\bm C}_{\rm ibest}$ and ${\bm \theta}_{\rm ibest}$ denote individual best variable sets of ${\bm C}$ and ${\bm \theta}$, respectively, and are of sizes of $(n+2)\times S$ and $m\times S$ as well; ${\bm C}_{\rm gbest} $ and ${\bm \theta}_{\rm gbest} $ denote the best variable set for ${\bm C}_{\rm ibest}$ and ${\bm \theta}_{\rm ibest}$, in other words, the global best variable set for ${\bm C}$ and ${\bm \theta}$, respectively. ${\bm C}_{v\min}$ and ${\bm C}_{v\max}$ denote the minimum and maximum searching velocities for ${\bm C}$, while ${\bm \theta}_{v\min}$ and ${\bm \theta}_{v\max}$ denote the minimum and maximum searching velocities for ${\bm \theta}$. Lines 12-22 in \textbf{Algorithm} \ref{alg.2} indicate that ${\bm C}_{\rm ibest} $ and ${\bm \theta}_{\rm ibest} $ always record individual best variable sets in computation history, ${\bm C}_{\rm gbest}$ and ${\bm \theta}_{\rm gbest}$ always record global best variable sets in the computation procedure. 

In \textbf{Algorithm} \ref{alg.3} and \ref{alg.5}, ${\bm C}_{v}$ and ${\bm \theta}_{v}$ denote the velocities for ${\bm C}$ and ${\bm \theta}$, respectively. \textbf{Algorithm} \ref{alg.4} indicate that ${\bm C}_{v}$ and ${\bm \theta}_{v}$ are updated using the velocity optimizers, and that ${\bm C}_{v}$ and ${\bm \theta}_{v}$, ${\bm C}$ and ${\bm \theta}$ are always constrained within the ranges bounded by the corresponding \textbf{Input} parameters. In \textbf{Algorithm} \ref{alg.3} and \ref{alg.5}, $r$, $r_{1}$, and $r_{2}$ denote three different random numerics within the range $(0,1)$ generated by computer.

In \textbf{Algorithm} \ref{alg.2}, $F_{N}({\bm C}_{\rm best},{\bm \theta}_{\rm gbest})$ and  $F_{N+1}({\bm C}_{\rm best},{\bm \theta}_{\rm gbest})$ denote the values of penalty function in Eq. (\ref{eq.8}) in $N^{\rm th}$ and $(N+1)^{\rm th}$ iteration, respectively. Lines 5 and 25 of \textbf{Algorithm} \ref{alg.2} indicate that the iteration would proceed if the difference of penalty function between two iterations is equal or greater than $VTR$ or the iteration rep is equal or smaller than the maximum iteration reps.

\textbf{Algorithm} \ref{alg.4} presents the linear updating strategies for the inertia, the individual learning factor, and the global learning factor, respectively. These three factors are applied to the modified particle swarm method in Lines 2 and 18 in \textbf{Algorithm} \ref{alg.5}. 

Lines 8-10 in \textbf{Algorithm} \ref{alg.2} compute the value of the penalty function with Eqs. (\ref{eq.6}) and (\ref{eq.8}), and consume most of the computation resources when excuting the pseudocodes. Therefore, the computation process to obtain the value of the penalty function is packed as a subroutine to facilitate parallel computation for time saving.

Lines 14 and 16 in \textbf{Algorithm} \ref{alg.3} illustrate two mutually exclusive initialization strategies. Line 14 is the variation of the initialization strategy recommended in Ref \cite{Zengguisen2019} ($\beta$ is a relatively small numeric), while Line 16 is a different initialization strategy with pure random numerics proposed by the authors. 

The pseudocodes in red in Line 19 of \textbf{Algorithm} \ref{alg.3} and Line 34 of \textbf{Algorithm} \ref{alg.5} sorts the elements in ${\bm \theta}[:,j]$ in the ascending order, and are the modification of the particle swarm method for the problem defined in this paper. Such modification is very slight comparing to the whole pseudocodes, but is very important to ensure solution convergence and global optimization. Though the last item in the penalty function ($\sum\limits_{i-1}^{m-1}\max(0,\theta_{i}-\theta_{i+1})^{2}$) in Eq. (\ref{eq.8}) theoretically considers the penalty components regarding to ${\bm \theta}$, Eq. (\ref{eq.8}) generally does not converge to the global minimum ($F({\bm X}^{*},{\bm \theta}^{*})=0$) in practical computation, if the the sortings are not included. In such a situation, the penalty function generally converges to some local minimum, and the value of $F({\bm X}^{*},{\bm \theta}^{*})$ is much greater than 0.

\subsection{Parameter selection}
\label{sec.4.2}

The range of ${\bm C}$ are chosen as:

\begin{equation}
  \label{eq.14}
  \left\{
    \begin{aligned}
      & {\bm C}_{\min}[1] = \frac{R_{\min}}{H_{\rm mid}+\sqrt{H_{\rm mid}^{2}-R_{\min}^{2}}} \\
      & {\bm C}_{\min}[2] = -2H_{\rm mid}\cdot\frac{1-{\bm C}_{\min}[1]^{2}}{1+{\bm C}_{\min}[1]^{2}} \\
      & {\bm C}_{\min}[i] = -{\rm e}^{2-i},\quad i=3,4,5,\cdots,n+2 \\
    \end{aligned}
  \right.
\end{equation}

\begin{equation}
  \label{eq.15}
  \left\{
    \begin{aligned}
      & {\bm C}_{\max}[1] = 1 \\
      & {\bm C}_{\max}[2] = 0 \\
      & {\bm C}_{\max}[i] = -{\bm C}_{\min}[i] ,\quad i=3,4,5,\cdots,n+2
    \end{aligned}
  \right.
\end{equation}

\noindent In Eq. (\ref{eq.14}), $H_{\rm mid}$ denotes the buried depth of the cavity, $R_{\min}$ denotes the radius of the incribed circle of the cavity, and these two parameters can be computed as

\begin{equation}
  \label{eq.16}
  \left\{
    \begin{aligned}
      & H_{\rm mid}=\frac{\left|y_{1}^{0}+y_{m}^{0}\right|}{2} \\
      & R_{\min} = \min\sqrt{{x_{i}^{0}}^{2}+\left(y_{i}^{0}-H_{\rm mid}\right)^{2}},\quad i=1,2,3,\cdots,m \\
    \end{aligned}
  \right.
\end{equation}

\noindent $R_{\min}$ is corresponding to the dashed black circle in Fig. \ref{fig.1}, thus, ${\bm C}_{\min}[1]$ and ${\bm C}_{\max}[1]$ always bound the cavity boundary between the black dashed circle and the ground surface. The average values of ${\bm C}_{\min}[i]$ and ${\bm C}_{\max}[i]$ ($i=2,3,4,\cdots,n+2$) are corresponding to those in Eq. (\ref{eq.9}). In Eq. (\ref{eq.14}), ${\rm e}$ denotes the natural logarithm base. The lower and upper limits for ${\bm C}[i] (i=3,4,5,\cdots,n+2)$ take the natural exponential form after many trial computations, however, such ranges are only for references, not the only option. Comparing to the initial values in Eq. (\ref{eq.9}), which are recommended in Ref \cite{Zengguisen2019}, the initial parameters in Eqs. (\ref{eq.14})-(\ref{eq.16}) are more excutable.

The range of ${\bm C}_{v\min}^{\rm init}$, ${\bm C}_{v\min}^{\rm end}$,${\bm C}_{v\max}^{\rm init}$, ${\bm C}_{v\max}^{\rm end}$, and ${\bm \theta}_{v\min}^{\rm init}$, ${\bm \theta}_{v\min}^{\rm end}$,${\bm \theta}_{v\max}^{\rm init}$, ${\bm \theta}_{v\max}^{\rm end}$ are chosen as:

\begin{equation}
  \label{eq.17}
  \left\{
    \begin{aligned}
      & {\bm C}_{v\min}^{\rm init}[i] =-0.02 \\
      & {\bm C}_{v\min}^{\rm end}[i] = -0.01 \\
      & {\bm C}_{v\max}^{\rm init}[i] =-{\bm C}_{v\min}^{\rm init}[i] \\
      & {\bm C}_{v\max}^{\rm end}[i] =-{\bm C}_{v\min}^{\rm end}[i] \\
    \end{aligned}
  \right.
  ,\quad
  i=1,2,3,\cdots,n+2
\end{equation}

\begin{equation}
  \label{eq.18}
  \left\{
    \begin{aligned}
      & {\bm \theta}_{v\min}^{\rm init}[i] = -0.01 \\
      & {\bm \theta}_{v\min}^{\rm end}[i] = -0.005 \\
      & {\bm \theta}_{v\max}^{\rm init}[i] = -{\bm \theta}_{v\min}^{\rm init}[i] \\
      & {\bm \theta}_{v\max}^{\rm end}[i] = -{\bm \theta}_{v\min}^{\rm end}[i] \\
    \end{aligned}
  \right.
  ,\quad
  i=1,2,3,\cdots,m
\end{equation}

\noindent Eqs. (\ref{eq.14}) and (\ref{eq.15}) indicate that the value ranges for ${\bm C}[i] (i=1,2,3,\cdots,n+2)$ and ${\bm \theta}[i] (i=1,2,3,\cdots,m)$ are very narrow, thus, the exploring velocities ${\bm C}_{v}$ and ${\bm \theta}_{v}$ should be relatively small, and should be further narrowed for algorithm stability when iteration proceeds.

The initial and ending values of the parameters of modified particle swarm method (including the inertia, individual learning, and global learning components) are chosen as follows:

\begin{equation}
  \label{eq.19}
  \left\{
    \begin{aligned}
      & w_{\rm init} = 0.7 \\
      & w_{\rm end} = 0.4 \\
    \end{aligned}
  \right.
  ,\quad
  \left\{
    \begin{aligned}
      & d_{1\rm init} = 2.5 \\
      & d_{1\rm end} = 0.5 \\
    \end{aligned}
  \right.
  ,\quad
  \left\{
    \begin{aligned}
      & d_{2\rm init} = 0.9 \\
      & d_{2\rm end} = 2.25 \\
    \end{aligned}
  \right.
\end{equation}

\noindent The other parameters are chosen as:

\begin{equation}
  \label{eq.20}
  \left\{
    \begin{aligned}
      & VTR = 10^{-6} \\
      & N_{\max} = 3000 \\
      & S = 10 \times (n+2+m) \\
      & \beta = 10^{-2} \\
    \end{aligned}
  \right.
\end{equation}

\subsection{Computation results and solution verification}
\label{sec.4.3}

Substituting the values in Eqs. (\ref{eq.14})-(\ref{eq.20}) and the four sets of initialization strategies into \textbf{Algorithm} \ref{alg.2} and encoding the pseudocodes with \texttt{FORTRAN} yields corresponding results. However, Lines 14 and 16 in \textbf{Algorithm} \ref{alg.3} indicate two mutually exclusive initial value strategies for ${\bm \theta}$. Meanwhile, whether or not to activate the pseudocodes marked in red to sort ${\bm \theta}[:,j]$ in Line 19 in \textbf{Algorithm} \ref{alg.3} and Line 34 in \textbf{Algorithm} \ref{alg.5} also results in two different sorting strategies for ${\bm \theta}[:,j]$. Therefore, four sets of initial value strategy combination can be established in Table \ref{tab.1}.

\begin{table}[htb]
  \centering
  \caption{Four sets of strategy combination}
  \begin{tabular}{ccc}
    \toprule
    Set & Initial value strategy of ${\bm \theta}$ in \textbf{Algorithm} \ref{alg.3} & Sorting or not? \\
    \midrule
    $A$ & Line 14  & Yes \\
    $B$ & Line 16  & Yes \\
    $C$ & Line 14  & No \\
    $D$ & Line 16  & No \\
    \bottomrule
  \end{tabular}
  \label{tab.1}
\end{table}

\begin{figure}[htb]
  \centering
  \begin{overpic}[abs,unit=1mm,grid=false]{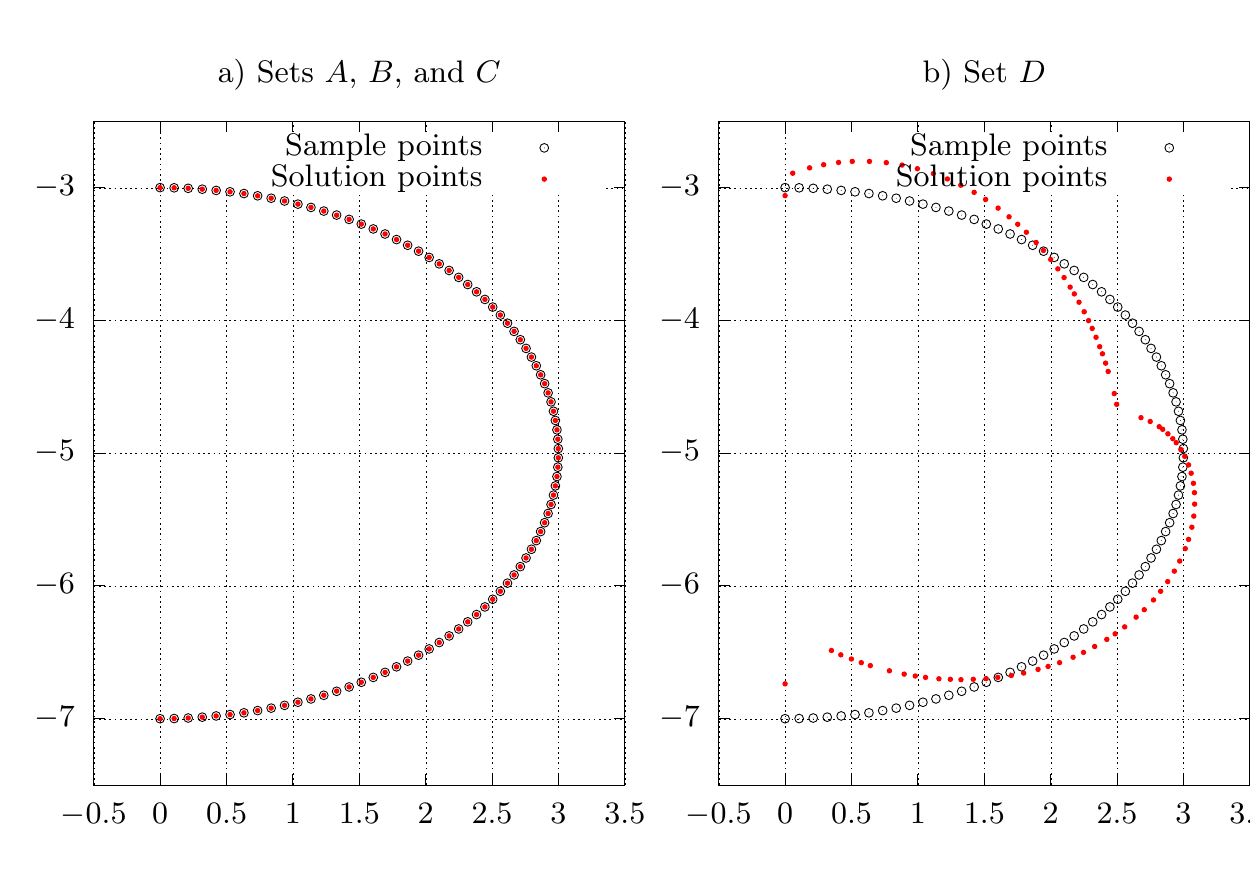}
  \end{overpic}
  \caption{Rectangular coordinate comparisons for the four sets}
  \label{fig.2}
\end{figure}

The computation results of rectangular coordinate comparison of the four sets can be seen in Fig. \ref{fig.2}. The coordinate comparisons of Sets $A$, $B$, and $C$ are almost the same, indicating that the previous three strategy combinations result in the same algorithm convergence. The coordinate comparison of Set $D$ is chaotic, indicating the initial value combination of Set $D$ is  unvaluable. The rectangular coordinate convergence procedure of Set $B$ can be seen in the \texttt{iteration.gif} file in the attachment.

\begin{figure}[htb]
  \centering
  \begin{overpic}[abs,unit=1mm,grid=false]{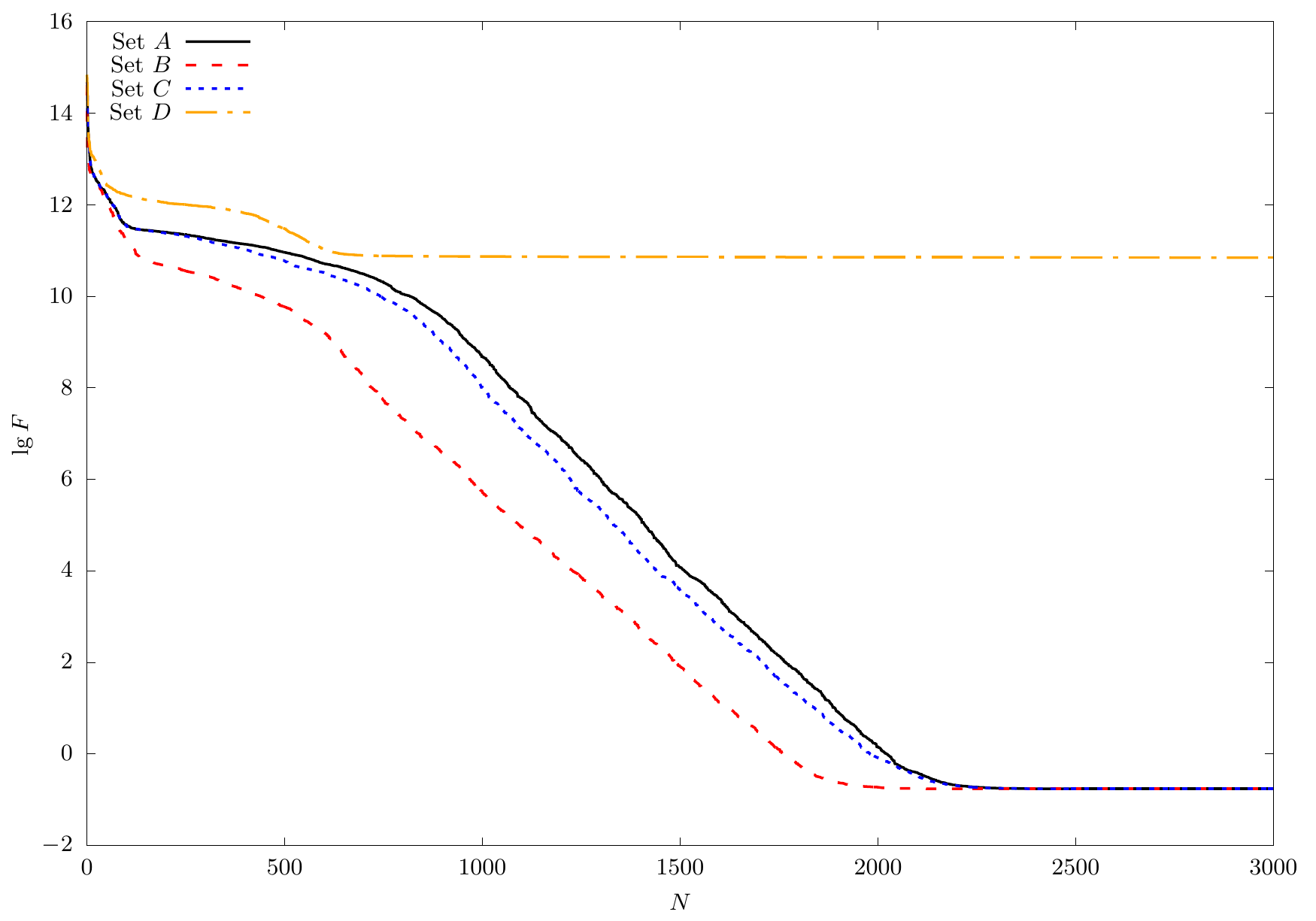}
    \put (20,80) {$ n = 10 $}
    \put (20,75) {$ m = 90 $}
    \put (20,70) {$ S = 1020 $}
    \put (20,65) {$ N_{\max} = 3000 $}
    \put (20,60) {$ VTR = 10^{-6} $}
    \put (20,55) {$ \beta = 0.01 $}
    \put (20,45) {$ H_{\rm mid} = 5 $}
    \put (20,40) {$ R_{\min} = 2 $}
    \put (20,30) {$ w_{\rm init} = 0.7,\quad w_{\rm end} = 0.4 $}
    \put (20,25) {$ d_{1\rm init} = 2.5,\quad d_{1\rm end} = 0.5 $}
    \put (20,20) {$ d_{2\rm init} = 0.9,\quad d_{2\rm end} = 2.25 $}

    \put (110,70) {$ {\bm C}_{v\min}^{\rm init} = -0.02,\quad {\bm C}_{v\min}^{\rm end} = -0.01 $}
    \put (110,60) {$ {\bm C}_{v\max}^{\rm init} = 0.02,\quad {\bm C}_{v\max}^{\rm end} = 0.01 $}
    \put (110,50) {$ {\bm \theta}_{v\min}^{\rm init} = -0.01,\quad {\bm \theta}_{v\min}^{\rm end} = -0.005 $}
    \put (110,40) {$ {\bm \theta}_{v\max}^{\rm init} = 0.01,\quad {\bm \theta}_{v\max}^{\rm end} = 0.005 $}
    \put (110,105) {
      \begin{tabular}{cc}
        Set & Elapsing time (s) \\
        \midrule
        $A$ & 27.904 \\
        $B$ & 26.478 \\
        $C$ & 18.826 \\
        $D$ & 19.206 \\
      \end{tabular}
    }
  \end{overpic}
  \caption{Penalty function value by logarithm ($\lg F$) against iteration reps ($N$)}
  \label{fig.3}
\end{figure}

The computation results of penalty function value by logarithm against iteration reps of all four sets can be seen in Fig. \ref{fig.3}. Clearly, the strategy combination of Set $D$ leads to some local minimum, just as mentioned in Section \ref{sec.4.1}, while the other three strategy combinations all lead to global minimum after about 2200 iteration reps. The consuming time of the convergence for each of these three gloabal minimums is acceptable (about 18-28 seconds to stop computation). The pseudocodes are encoded via \texttt{FORTRAN} using \texttt{NVFORTRAN} compiler with assistance of \texttt{OPENACC} and \texttt{CUDA-11.7} accelaration frameworks, and are performed on a workstation with a CUDA GPU (GeForce RTX 3060 Lite Hash Rate).

Fig. \ref{fig.3} indicates that Set $B$ would reach convergence with fewer iteration reps than Sets $A$ and $C$. One possible reason is that the randomness is included in the initial values to slightly decrease the numerical differences between the initial values and the final values, when comparing to the uniform initial value strategy of Sets $A$ and $C$. Meanwhile, the sorting strategy further ensures such an advantage is kept until the optimization reaches convergence.
 
Till now, the solution method and procedure of the optimization of the penalty function in Eq. (\ref{eq.8}) are presented in detail, and primary aim of this paper has been achieved.

\subsection{Solution accuracy and efficiency}
\label{sec.4.4}

The optimization procedure for the global minimum of Eq. (\ref{eq.8}) is actually a fitting procedure of sequential points. For a fitting procedure, both underfitting and overfitting should be avoided. To be specific, the combination $n$ and $m$ is dominant on whether or not the fitting procedure of Eq. (\ref{eq.8}) would fall into underfitting or overfitting, as long as the other parameters remain the same. Furthermore, the combination of $n$ and $m$ also determines the solution accuracy and efficieny of Eq. (\ref{eq.8}). Therefore, we conduct a series of numerical experiments to determine the relatively reasonable combination of $n$ and $m$. 

In the experiments, both $n$ and $m$ monotonously and discretely increase, and finally 12 sets of $n$ and $m$ combinations are selected, as illustrated in Table \ref{tab.2}. When substituting the 12 sets of $n$ and $m$ into the solution, the maximum iteration rep $N_{\max}$ slowly increases by an interval of 100 from 0, until a stable value of $\lg[F_{N+1}({\bm C}_{\rm gbest},{\bm \theta}_{\rm gbest})]$ is found. The solution results of all 12 sets of $n$ and $m$ combinations are shown in Table \ref{tab.2}.

The third and fifth columns in Table \ref{tab.2} show that as the combination of $n$ and $m$ increases monotonously, $N_{\max}$ and elapsing time both generally increase as well. The fourth column in Table \ref{tab.2} indicates that there exists a narrow range to obtain global minimum for combination of $n$ and $m$ (Sets 9 and 10), within which both accuracy and efficiency would be reansonable. When the combination of $n$ and $m$ is out of such a solution range, local minimum would be found due to underfitting or overfitting, instead of desired global minimum. Table \ref{tab.2} further indicates that though global minimum for Eq. (\ref{eq.8}) theoretically exists, tangibly locating and finding out that global minimum may be time-consuming.

\begin{table}[hbt]
  \centering
  \caption{Computation results for different combinations of $n$ and $m$}
  \begin{tabular}{cccccc}
    \toprule
    Set & $n$ & $m$ & $N_{\max} $ & $\lg F$ (stable)  & Elapsing time (s) \\
    \midrule
    1 & 1 & 3 & 100 & 8.803 & 0.107 \\
    2 & 1 & 6 & 300 & 8.686 & 0.147 \\
    3 & 1 & 10 & 400 & 8.935 & 0.271 \\
    4 & 2 & 10 & 400 & 7.747 & 0.271 \\
    5 & 3 & 10 & 500 & 6.562 & 0.434 \\
    6 & 5 & 15 & 1100 & 4.268 & 1.110 \\
    7 & 5 & 30 & 1100 & 4.671 & 1.460 \\
    8 & 10 & 30 & 2000 & 9.065 & 36.619 \\
    9 & 10 & 60 & 2900 & -0.934 & 12.358 \\
    10 & 10 & 90 & 2200 & -0.739 & 19.188 \\
    11 & 15 & 90 & 4000 & 10.214 & 43.418 \\
    12 & 15 & 180 & 5000 & 9.500 & 182.543 \\
    \bottomrule                       
  \end{tabular}
  \label{tab.2}
\end{table}

\section{Solution comments and further discussion}
\label{sec.5}

In this paper, the following works have been done:

\noindent (1) The penalty function for the solution of the nonlinear equation system in Eq. (\ref{eq.6}) are rephrased with detailed mathematical discription.

\noindent (2) The nonconvexity of the penalty function has been proven via contradiction proof based on the definition of convex function, and detailed pseudocodes and corresponding \texttt{FORTRAN} codes are provided.

\noindent (3) Due to the nonconvexity of the penalty function, the modified Particle Swarm Method has been presented very detailedly for the solution to ensure that the solution procedure can be fully repeated and possibly further improved, and the pseudocodes and corresponding \texttt{FORTRAN} codes have been presented as well. 

\noindent (4) The initial values of the coefficients of the conformal mapping function (Eqs. (\ref{eq.14})-(\ref{eq.16})) have been modified from single values in Ref \cite{Zengguisen2019} to value ranges to ensure that the final results of these coefficients are always within the value ranges to reach convergence.

\noindent (5) Four strategies for the initial values of the polar angles of the sample points have been compared, and a new initial value strategy has been found, which is more efficient than the one recommended in Ref \cite{Zengguisen2019}.

\noindent (6) Solution for global minimum can only be found within solution range to avoid underfitting or overfitting with satisfying accuracy and efficiency.

Though the conformal mapping function in Ref \cite{Zengguisen2019} greatly extends the understanding and usage of the complex variable method, such mapping still has its limits. As can be seen in both Ref \cite{Zengguisen2019} and this paper, the solution procedure only focuses on the right half of the cavity in the physical plane, since the conformal mapping function is only suitable to symmetrical cavities in lower-half plane, and the axis of symmetry is $x=0$ (Fig. \ref{fig.1}a). Such property can be proven by simply substituting $-\theta$ into Eq. (\ref{eq.3}), since the range of $\theta$ in $[0,2\pi]$ can be rephased as $[-\pi,\pi]$ due to periodicity of trigonometric function:

\begin{equation}
  \left\{
    \begin{aligned}
      & x(\rho,-\theta)=-c_{0} \frac{2\rho\sin(-\theta)}{1+\rho^2-2\rho\cos(-\theta)}-\sum\limits_{k=1}^{n}{c_{k}}\left(\rho^k+\rho^{-k}\right)\sin{k}(-\theta) = - x(\rho,\theta) \\
      & y(\rho,-\theta)=c_{0} \frac{1-\rho^2}{1+\rho^2-2\rho\cos(-\theta)}+\sum\limits_{k=1}^{n}{c_{k}}\left(\rho^k-\rho^{-k}\right)\cos{k}(-\theta) = y(\rho,\theta)
    \end{aligned}
  \right.
  \label{eq.21}
\end{equation}

Eq. (\ref{eq.21}) indicates that $\theta=0$ always results in $x=0$ for arbitrary $\rho$. With no doubts, the symmetry of the conformal mapping function would simplify the solution procedure to obtain the unknown coefficients. However, on the other hand, symmetrical cavities caused by excavation are rarely seen in real tunnel enginnering, since over excavation generally exists, and the symmetry prevents further usage of the conformal mapping function and complex variable method. Therefore, an interesting and valuable future research direction is to propose a new conformal mapping function for asysmetrical cavities in a lower-half plane, which would certainly further extend the usage of the complex variable method. 

\clearpage
\section*{Acknowlegement}
\label{acknowlegement}

This study is financially supported by Scientific Research Foundation of Fujian University of Technology (Grant No. GY-Z20094), the National Natural Science Foundation of China (Grant No. 52178318), and Education Foundation of Fujian Province (Grant No. JAT210287). The authors would like to thank Professor Changjie Zheng, Ph.D. Yiqun Huang, and Ph.D. Xiaoyi Zhang for their suggestions on this study.

\clearpage
\begin{algorithm}
  \caption{Verification of nonconvexity for Eq. (\ref{eq.8})}
  \label{alg.1}
  \begin{algorithmic}[1]
    \REQUIRE$A$, $B$, $H$, $n$, $m$
    \ENSURE$N_{\max} $, $G_{\rm minus} $
    \STATE Computing $x^{0}_{i}$, $y^{0}_{i}$($i=1,2,3,\cdots,m$) in Eq. (\ref{eq.11})
    \STATE$N \leftarrow 1$
    $\quad$//$N$ is an integer, which denotes the iteration reps
    \STATE$G \leftarrow 1.0$
    
    \WHILE{$G \geq 0$}
      \FOR{each $j\in[1,3]$}
        \STATE${\bm C}[:,j] \leftarrow 0$, ${\bm \theta}[:,j] \leftarrow 0$
          $\quad$//Always initializing ${\bm C}$, ${\bm \theta}$
        \STATE$F({\bm C}[:,j],{\bm \theta}[:,j]) \leftarrow 0$
          $\quad$//Always initializing $F({\bm C},{\bm \theta})$ in Eq. (\ref{eq.8})
      \ENDFOR
        $\quad$//Initialization finished
    
      \FOR{each $j\in[1,2]$}
        \STATE${\bm C}[1,j] \leftarrow t$
        \STATE${\bm C}[2,j] \leftarrow (2t-1) \cdot 20$
          \FOR{each $i\in[3,n+2]$}
            \STATE${\bm C}[i,j] \leftarrow 2t-1$
          \ENDFOR
          \STATE{${\bm \theta}[1,j] \leftarrow 0$}
          \FOR{each $i\in[2,m-1]$}
            \STATE${\bm \theta}[i,j] \leftarrow t\cdot\pi$
          \ENDFOR
          \STATE{${\bm \theta}[m,j] \leftarrow \pi$}
        \STATE Sorting ${\bm \theta}[:,j]$ in the ascending order
      \ENDFOR
        $\quad$//Asignment of ${\bm C}[:,1] $ and ${\bm C}[:,2] $, ${\bm \theta}[:,1] $ and ${\bm \theta}[:,2] $ finished
    
      \STATE${\bm C}[:,3] \leftarrow \lambda\cdot{\bm C}[:,1]+(1-\lambda)\cdot{\bm C}[:,2]$
      \STATE{${\bm \theta}[:,3] \leftarrow \lambda \cdot {\bm \theta}[:,1] + (1-\lambda) \cdot {\bm \theta}[:,2] $}
        $\quad$//Generating the mediate variable in Eq. (\ref{eq.9})
      \STATE Computing $F({\bm C}[:,1],{\bm \theta}[:,1])$, $F({\bm C}[:,2],{\bm \theta}[:,2])$, $F({\bm C}[:,3],{\bm \theta}[:,3])$ via Eq. (\ref{eq.8})
      \STATE Computing $G$ via Eq. (\ref{eq.9})
      \STATE$N \leftarrow N+1$
        $\quad$//updating $N$
      \PRINT $N$, $G$
    \ENDWHILE
    \STATE{$N_{\max} \leftarrow N$, $G_{\rm minus} \leftarrow G$}
    \PRINT{$N_{\max}$, $G_{\rm minus}$}
  \end{algorithmic}
\end{algorithm}

\clearpage
\begin{algorithm}
  \caption{Modified particle swarm method}
  \label{alg.2}
  \begin{algorithmic}[1]
    \REQUIRE$VTR$, $N_{\max}$, $S$, $\varOmega$, $\alpha$, $c_{0}$,
      ${{\bm C}_{v\min}^{\rm init}}$, ${{\bm C}_{v\min}^{\rm end}}$, ${{\bm C}_{v\max}^{\rm init}}$, ${{\bm C}_{v\max}^{\rm end}}$,
      ${{\bm \theta}_{v\min}^{\rm init}}$, ${{\bm \theta}_{v\min}^{\rm end}}$, ${{\bm \theta}_{v\max}^{\rm init}}$, ${{\bm \theta}_{v\max}^{\rm end}}$,
      $w_{\rm init}$, $w_{\rm end}$,
      $d_{1 \rm init}$, $d_{1 \rm end}$, $d_{2 \rm init}$, $d_{2 \rm end}$
      
    \ENSURE ${\bm C}^{*}$, ${\bm \theta}^{*} $
    \STATE{\textbf{CALL} \textbf{Algorithm} \ref{alg.3}
      $\quad$//Variable initialization}
    \STATE{$N \leftarrow 0$} $\quad$//Initializing iteration number
    \STATE$F_{N} ({\bm C}_{\rm gbest} , {\bm \theta}_{\rm gbest})  \leftarrow \varOmega^{3} $
    \STATE$F_{N+1} ({\bm C}_{\rm gbest} , {\bm \theta}_{\rm gbest}) \leftarrow \varOmega^{2} $
    \WHILE{$F_{N}({\bm C}_{\rm gbest} ,{\bm \theta}_{\rm gbest})-F_{N+1}({\bm C}_{\rm gbest} , {\bm \theta}_{\rm gbest})\geq VTR$ \OR $N \leq N_{\max}$}
      \STATE{\textbf{CALL} \textbf{Algorithm} \ref{alg.4}}
      \STATE{$F_{N} ({\bm C}_{\rm gbest} ,{\bm \theta}_{\rm gbest}) \leftarrow F_{N+1}({\bm C}_{\rm gbest} ,{\bm \theta}_{\rm gbest})$}
       \FOR{each $j\in[1,S]$}
        \STATE Computing $F({\bm C}[:,j], {\bm \theta}[:,j])$ and
          $F({\bm C}_{\rm ibest}[:,j],{\bm \theta}_{\rm ibest}[:,j])$
      \ENDFOR $\quad$//\emph{Parallel computation}
      \FOR{each $j\in[1,S]$}
        \IF{$F({\bm C}[:,j],{\bm \theta}[:,j]) <
          F({\bm C}_{\rm ibest} [:,j],{\bm \theta}_{\rm ibest} [:,j])$}
          \STATE${\bm C}_{\rm ibest}[:,j] \leftarrow {\bm C}[:,j]$
          \STATE${\bm \theta}_{\rm ibest}[:,j] \leftarrow {\bm \theta}[:,j]$
        \ENDIF
        \IF{$F({\bm C}[:,j] ,{\bm \theta}[:,j]) < F({\bm C}_{\rm gbest},{\bm \theta}_{\rm gbest})$}
          \STATE ${\bm C}_{\rm gbest} \leftarrow {\bm C}_{\rm ibest}[:,j] $
          \STATE${\bm \theta}_{\rm gbest} \leftarrow {\bm \theta}_{\rm ibest}[:,j] $
          \STATE$F_{2} ({\bm C}_{\rm gbest},{\bm \theta}_{\rm gbest})
              \leftarrow F({\bm C}_{\rm ibest} [:,j],{\bm \theta}_{\rm ibest} [:,j])$
        \ENDIF
        \STATE{\textbf{CALL} \textbf{Algorithm} \ref{alg.5}
          $\quad$//Updating and constraining ${\bm C}_{v}$ and ${\bm C}$, ${\bm \theta}_{v}$ and ${\bm \theta}$}
      \ENDFOR
      \PRINT{$N$, $\lg\left[F_{N+1} ({\bm C}_{\rm gbest},{\bm \theta}_{\rm gbest})\right]$}
      \STATE{$N \leftarrow N+1$}
    \ENDWHILE
    \STATE{${\bm C}^{*} \leftarrow {\bm C}_{\rm gbest}$}
    \STATE{${\bm \theta}^{*} \leftarrow {\bm \theta}_{\rm gbest}$}
    \PRINT{$ {\bm C}^{*}$, ${\bm \theta}^{*}$}
  \end{algorithmic}
\end{algorithm}

\clearpage
\begin{algorithm}
  \caption{Modified particle swarm method: variable initialization}
  \label{alg.3}
  \begin{algorithmic}[1]
    \FOR{each $j\in[1,S]$}
      \FOR{each $i\in[1,n+2]$}
        \STATE${\bm C}_{v}[i,j] \leftarrow 0$
      \ENDFOR
      \FOR{each $i\in[1,m]$}
        \STATE${\bm \theta}_{v}[i,j] \leftarrow 0$
      \ENDFOR
    \ENDFOR
    \FOR{each $j\in[1,S]$}
      \FOR{each $i\in[1,n+2]$}
        \STATE${\bm C}[i,j] \leftarrow {\bm C}_{\min}[i]+
          \left({\bm C}_{\max}[i] -{\bm C}_{\min}[i] \right)\cdot r$
      \ENDFOR
      \FOR{each $i\in[1,m]$}
        \STATE${\bm \theta}[i,j] \leftarrow \pi\cdot\frac{i-1}{m-1}+\pi\cdot \frac{\beta}{m-1} \cdot (2r-1)$ 
        \STATE\OR
        \STATE{${\bm \theta}[i,j] \leftarrow \pi\cdot r$}
      \ENDFOR
        $\quad$//Choosing either strategy for initialization 
      \STATE{${\bm \theta}[1,j] \leftarrow 0$, ${\bm \theta}[m,j] \leftarrow \pi$}
      \STATE \textcolor{red}{Sorting ${\bm \theta}[:,j]$ in the ascending order}
        $\quad$//Sorting to guarantee convergence and to accelarate the algorithm
    \ENDFOR
    \STATE${\bm C}_{\rm ibest} \leftarrow {\bm C}$, ${\bm \theta}_{\rm ibest} \leftarrow {\bm \theta}$
    \FOR{each $i\in[1,n+2]$}
      \STATE${\bm C}_{\rm gbest}[i] \leftarrow \varOmega$
    \ENDFOR
    \FOR{each $i\in[1,m]$}
      \STATE${\bm \theta}_{\rm gbest}[i] \leftarrow \varOmega$
    \ENDFOR
    \FOR{each $j\in[1,S]$}
      \STATE$F({\bm C}[:,j],{\bm \theta}[:,j]) \leftarrow \varOmega$
      \STATE$F({\bm C}_{\rm ibest}[:,j] ,{\bm \theta}_{\rm ibest}[:,j] ) \leftarrow \varOmega$
    \ENDFOR
  \end{algorithmic}
\end{algorithm}

\clearpage
\begin{algorithm}
  \caption{Modified particle swarm method: updating the parameters related to iteration reps}
  \label{alg.4}
  \begin{algorithmic}[1]
    \STATE$w \leftarrow w_{\rm init}+(w_{\rm end}-w_{\rm init})\cdot\frac{N}{N_{\max}} $
    \STATE$d_{1}\leftarrow d_{1\rm init}+(d_{1\rm end}-d_{1\rm init})\cdot\frac{N}{N_{\max}} $
    \STATE$d_{2}\leftarrow d_{2\rm init}+(d_{2\rm end}-d_{2\rm init})\cdot\frac{N}{N_{\max}} $
    \FOR{each $i\in[1,n+2]$}
    \STATE{${\bm C}_{v\min}[i]={\bm C}_{v\min}^{\rm init}[i]+({\bm C}_{v\min}^{\rm end}[i]-{\bm C}_{v\min}^{\rm init}[i])\cdot\frac{N}{N_{\max}} $}
    \STATE{${\bm C}_{v\max}[i]={\bm C}_{v\max}^{\rm init}[i]+({\bm C}_{v\max}^{\rm end}[i]-{\bm C}_{v\max}^{\rm init}[i])\cdot\frac{N}{N_{\max}} $}
    \ENDFOR
    \FOR{each $i\in[1,m]$}
    \STATE{${\bm \theta}_{v\min}[i]={\bm \theta}_{v\min}^{\rm init}[i]+({\bm \theta}_{v\min}^{\rm end}[i]-{\bm \theta}_{v\min}^{\rm init}[i])\cdot\frac{N}{N_{\max}} $}
    \STATE{${\bm \theta}_{v\max}[i]={\bm \theta}_{v\max}^{\rm init}[i]+({\bm \theta}_{v\max}^{\rm end}[i]-{\bm \theta}_{v\max}^{\rm init}[i])\cdot\frac{N}{N_{\max}} $}
    \ENDFOR
  \end{algorithmic}
\end{algorithm}

\clearpage
\begin{algorithm}
  \caption{Modified particle swarm method: updating and constraining ${\bm C}_{v}$ and ${\bm C}$, ${\bm \theta}_{v}$ and ${\bm \theta}$}
  \label{alg.5}
  \begin{algorithmic}[1]
    \FOR{each $i\in[1,n+2]$}
      \STATE{${\bm C}_{v} [i,j] \leftarrow
        w \cdot {\bm C}_{v}[i,j]
        + d_{1} \cdot r_{1} \cdot ({\bm C}_{\rm ibest}[i,j] - {\bm C}[i,j])
        + d_{2} \cdot r_{2} \cdot ({\bm C}_{\rm gbest}[i] - {\bm C}[i,j])$}
      \IF{${\bm C}_{v}[i,j]>{\bm C}_{v\max}[i]$}
        \STATE{${\bm C}_{v}[i,j]\leftarrow{\bm C}_{v\max}[i]$}
      \ELSIF{${\bm C}_{v}[i,j]<{\bm C}_{v\min}[i]$}
        \STATE{${\bm C}_{v}[i,j]\leftarrow{\bm C}_{v\min}[i]$}
      \ENDIF
    \ENDFOR
      $\quad$//Updating and constraining ${\bm C}_{v} $
    \FOR{each $i\in[1,n+2]$}
      \STATE{${\bm C}[i,j] \leftarrow {\bm C}[i,j] + {\bm C}_{v}[i,j] $}
      \IF{${\bm C}[i,j]>{\bm C}_{\max}[i]$}
        \STATE{${\bm C}[i,j]\leftarrow{\bm C}_{\max}[i]$}
      \ELSIF{${\bm C}[i,j]<{\bm C}_{\min}[i]$}
        \STATE{${\bm C}[i,j]\leftarrow{\bm C}_{\min}[i]$}
      \ENDIF
    \ENDFOR
    $\quad$//Updating and constraining ${\bm C}$
    \FOR{each $i\in[1,m]$}
      \STATE{${\bm \theta}_{v}[i,j] \leftarrow
        w \cdot {\bm \theta}_{v}[i,j]
        + d_{1}\cdot r_{1}\cdot({\bm \theta}_{\rm ibest}[i,j]-{\bm \theta}_{\rm ibest}[i,j])
        + d_{2}\cdot r_{2}\cdot({\bm \theta}_{\rm gbest}[i]-{\bm \theta}_{\rm ibest}[i,j])$}
      \IF{${\bm \theta}_{v}[i,j]>{\bm \theta}_{v\max}[i]$}
        \STATE{${\bm \theta}_{v}[i,j]\leftarrow{\bm \theta}_{v\max}[i]$}
      \ELSIF{${\bm \theta}_{v}[i,j]<{\bm \theta}_{v\min}[i]$}
        \STATE{${\bm \theta}_{v}[i,j]\leftarrow{\bm \theta}_{v\min}[i]$}
      \ENDIF
    \ENDFOR
      $\quad$//Updating and constraining ${\bm \theta}_{v}$
    \FOR{each $i\in[1,m]$}
      \STATE{${\bm \theta}[i,j]\leftarrow{\bm \theta}[i,j]+{\bm \theta}_{v}[i,j]$}
      \IF{${\bm \theta}[i,j]>{\bm \theta}_{\max}[i]$}
        \STATE{${\bm \theta}[i,j]\leftarrow{\bm \theta}_{\max}[i]$}
      \ELSIF{${\bm \theta}[i,j]<{\bm \theta}_{\min}[i]$}
        \STATE{${\bm \theta}[i,j]\leftarrow{\bm \theta}_{\min}[i]$}
      \ENDIF
    \ENDFOR
      $\quad$//Updating and constraining ${\bm \theta}$
    \STATE{${\bm \theta}[1,j] \leftarrow 0$, ${\bm \theta}[m,j] \leftarrow \pi$}
      $\quad$//Mandatorily constrain upper and lower limits of ${\bm \theta}[:,j]$
    \STATE{\textcolor{red}{Sorting ${\bm \theta}[:,j]$ in the ascending order}}
      $\quad$//Sorting to ensure convergence and to accelarate the algorithm
  \end{algorithmic}
\end{algorithm}

\clearpage

\clearpage
 
\clearpage

\clearpage
\bibliographystyle{plain} 
\bibliography{paper3.bib}

\end{document}